\documentclass{amsart}

\usepackage{epsfig,amsmath,graphics}
\usepackage[psamsfonts]{amssymb}
\usepackage[latin1]{inputenc}
\usepackage[usenames]{color}

\newtheorem{definition}{\bf Definition}
\newtheorem{theorem}{\bf Theorem}
\newtheorem{lemma}{\bf Lemma}

\renewcommand{\epsilon}{\varepsilon} 
\renewcommand{\emptyset}{\varnothing}

\def\tend{\longrightarrow}

\def\ds{\displaystyle}

\def\on{\operatorname}

\def\C{{\mathbb C}}
\def\D{{\mathbb D}}

\def\N{{\mathbb N}}

\def\R{{\mathbb R}}
\def\T{{\mathbb T}}
\def\U{{\mathbb U}}
\def\Z{{\mathbb Z}}

\def\cal{\mathcal}

\def\proof{\noindent{\bf Proof. }}
\def\remark{\vskip.2cm \noindent{\bf Remark. }}

\begin{document}

\title{Quadratic Siegel Disks with Rough Boundaries.}
\begin{author}[X.~Buff]{Xavier Buff}
\email{buff$@$picard.ups-tlse.fr}
\address{ %
  Universit\'e Paul Sabatier\\ 
  Laboratoire Emile Picard \\
  118, route de Narbonne \\
  31062 Toulouse Cedex \\
  France }
\end{author}
\begin{author}[A.~Ch\'eritat]{Arnaud Ch\'eritat}
\email{cheritat$@$picard.ups-tlse.fr}
\address{ %
  Universit\'e Paul Sabatier\\ 
  Laboratoire Emile Picard \\
  118, route de Narbonne \\
  31062 Toulouse Cedex \\
  France }
\end{author}

\begin{abstract}
In the quadratic family (the set of polynomials of degree $2$),
Petersen and Zakeri \cite{PZ} proved the existence of Siegel disks
whose boundaries are Jordan curves, but not quasicircles. In their
examples, the critical point is contained in the curve.

In the first part, as an illustration of the flexibility of the tools
developped in~\cite{BC1}, we prove the existence of examples that do not
contain the critical point.

In the second part, using a more abstract point of view
(suggested by Avila in~\cite{A}), we show that we can control quite
precisely the degree of regularity of the boundary of the Siegel disks
we create by perturbations.

\end{abstract}

\maketitle

Notations: $\U$ is the set of complex numbers with norm $=1$. Here, $r
\U$ will be used as a shorthand for: the circle of center $0$ and
radius $r$. Sometimes, we will note $(r,t)$ for $r \exp(i 2\pi
t)$. The symbol $\T$ will denote the quotient $\R/\Z$. For $\theta \in
\R$, 
\[P_\theta(z) = \exp(i2\pi \theta) z +z^2.\]
If the fixed point $z=0$ is
linearizable, then $\Delta(\theta)$ is
the Siegel disk of $P_\theta$ at $0$, and $r(\theta)$ its conformal
radius. Otherwise $\Delta(\theta) = \emptyset$ and $r(\theta) = 0$.

Reminder: $\forall \theta\in\R$, $r(\theta)<4$.

\

Let $\cal D_2$ denote the set of bounded type irrational numbers.
Let us recall the following

\begin{lemma}[Herman]\label{lem_Herman}
For all $\theta\in\cal D_2$, assume that $U$ is a
connected open set containing $0$ and $f : U \to \C$ is a
holomorphic function which fixes $0$ with derivative $e^{2\pi i
\theta}$. Let $\Delta$ be the Siegel disk of $f$ at $0$ (which exists by a
theorem of Siegel). If $U$ is simply connected and $f$ is univalent,
then $\Delta$ cannot have compact closure in $U$.
\end{lemma}

\remark In fact, Herman's theorem is stronger: $U$ needs not to be
simply connected, and the condition on $\theta$ is weaker (it is
called the Herman condition). See~\cite{He} for a reference.

\

An essential tool is the following
\begin{lemma}[Buff, Chéritat] (independently by A.~Avila)\label{lem_abc}
  For all $\theta$ Bruno and $r < r(\theta)$, there exists a
  sequence of $\cal D_2$ numbers $\theta_n \tend \theta$ such
  that $r(\theta_n) \tend r$.
\end{lemma}
Both \cite{BC1} and \cite{A} have stronger statements.
Note that we do not only require $\theta_n$ to be Bruno, but also to be
a bounded type number (the bound varies with $n$). To the
interested reader, we recommend \cite{A} for its simplicity (the fact
that $\theta_n$ can be taken $\cal D_2$ is not explicitly stated, but it
follows rather easily).
However, \cite{BC1} also provides a \emph{small cycle}.

\part{Rough Siegel disks}

\begin{definition}
  Let $L\subset \C$ be a Jordan curve. For any pair $x$,$y\in L$ of
  distinct points, let $U$ and $V$ be the connected components of
  $L-\{x,y\}$ and define the \emph{pinching} by
  \[\on{pinch}(L,x,y) = \min\big(\on{diam}(U),\on{diam}(V)\big) /
  \on{dist}(x,y).\]
  A \emph{quasicircle} is a Jordan curve whose pinchings are bounded over all
  possible pairs $x,y$.
  If the bound is $K$, we will say that we have a $K$-quasicircle.
\end{definition}

For instance a round circle is a $1$-quasicircle.
Note that a $1$-quasicicle is not necessarily a round circle.

\begin{lemma}\label{lem_geom_1}
Assume that $K>1$, and that $L_n$ are
$K$-quasicircles. Let $U_n$ be the bounded component of $\C-L_n$. Assume that
the limit of $U_n$ for the Carathéodory topology on domains containing
$0$ is equal to $\D$.
Then $L_n$ has Hausdorff limit equal to $\U$.
\end{lemma}
We leave it as an exercise to the reader. The assumption just means
that every compact subset of $\D$ is eventually contained in $U_n$,
and that every point of $\U$ has distance to $\C-U_n$ tending to $0$.

\begin{lemma}\label{lem_geom_2}
Assume that $K>1$, and that $L_n$ are $K$-quasicircles. Let $U_n$ be
the bounded component of $\C-L_n$. Assume that $0\in U_n$ and that
$L_n$ has Hausdorff 
limit equal to $\U$. Then the conformal mapping $f_n$ from $\D$
to $U_n$ mapping $0$ to $0$ with real positive derivative tends uniformly on
$\D$ to the identity.
\end{lemma}
\proof
A well known property
states that there exists $K'$
depending only on $K$ such that $f_n$ extends to a $K'$-quasiconformal
homeomorphism of $\C$. The set of $K'$-quasiconformal homeomorphisms
of $\C$ fixing $0$ and $1$ is compact. Since $f_n(1)$ tends to $\U$,
$f_n$ lies in a compact set. Any limit value of the
sequence $f_n$ must map $\D$ to itself, fixing $0$ with real positive
derivative, and thus is the identity on $\overline{\D}$.
\qed

\begin{lemma}\label{lem_C1}
Assume that $f,f_n : \T\to \C$ are $C^1$ functions, that the derivative of
$f$ does not vanish, that $f$ is injective, and that $f(\T)$ is a
$K$-quasicircle. Assume that the sequence $f_n-f$ and its derivative
tend uniformly to $0$. Then $\forall K'>K$, $\exists N\in\N$ $\forall
n \geq N$, $f_n(\T)$ is a Jordan curve and a $K'$-quasicircle.
\end{lemma}
\proof
  Let $m=\min |f'|$ and $m_n = \min |f'_n|$. Then $m>0$ and $m_n>0$.
  Let $\epsilon'>0$ such that $d(t,t')<\epsilon' \Longrightarrow
  |f'(t)-f'(t')|< m/8$. Let $N$ big enough so that $\forall n \geq N$,
  $|f'_n - f'| < m/8$, hence $\frac{7}{8} m <  m_n < \frac{9}{8}
  m$. Then $\forall n 
  \geq N$ and $d(t,t')<\epsilon$, 
  we have $|f'_n(t)-f'_n(t')| < m_n /2$. This implies
  $t,t'$ cannot parameterize a pinching of $f_n(\T)$ greater
  than $1$ if $d(t,t')<\epsilon$.
  Now, if we let $\mu = \min|f(t)-f(t')|$ taken on all pairs $t,t' \in
  \T$ with  $d(t,t') \geq \epsilon$, we have
  $|f_n -f| < \eta < \mu/2 \Longrightarrow f_n(\T)$ is a
  $K'$-quasicircle with \[K' = K \frac{\mu+2\eta}{\mu-2\eta}.\]
\qed

\

The key lemma is the following.
\begin{lemma}[perturbation]\label{lem_maillon}
  For all $\theta$ Bruno, for all $r_1<r_2<r(\theta)$, for all $K>1$ and for
  all $\epsilon>0$ there exists a Bruno number $\theta'$ and $r'>0$ such that
  \begin{enumerate}
   \item $|\theta'-\theta|<\epsilon$
   \item $r_1<r'<r_2$
   \item $r(\theta') > r'$
   \item $\| \phi_{\theta} - \phi_{\theta'} \|_\infty < \epsilon$ on
   the circle $r'\U$ \label{item_inf}
   \item the image of this circle by $\phi_{\theta'}$
   has a pinching $>K$
  \end{enumerate}
\end{lemma}
\proof
Let $r_3 = \frac{r_1+r_2}{2}$.
If $\phi_{\theta}(r' \U)$ has a pinching $> K$ for some $r'\in ]r_1,r_3[$,
then we are done with $\theta'= \theta$. Otherwise, let $\theta_n
\tend \theta$ provided by lemma~\ref{lem_abc} such that 
$r(\theta_n) \tend r_3$.
We claim that for all $n$ big enough, there exists an $r\in
[r_1,r(\theta_n)[$ such that 
$\phi_{\theta_n}(r\U)$ has a pinching $\geq K+1$.
Otherwise, lemma~\ref{lem_geom_1} would imply that
$\partial \Delta(\theta_n)$ tends to $\phi_\theta(r_3 \U)$ and thus is
eventually contained in $\Delta(\theta)$, contradicting
lemma~\ref{lem_Herman} since $\theta_n \in \cal D_2$.

Let then $r'_n$ be the infimum of the set of $r \in [r_1,r(\theta_n)[$ such
that $\phi_{\theta_n} (r\U)$ has a pinching $\geq K+1$.
Then $\phi_{\theta_n}$ has a $(K+1)$-pinching, otherwise
lemma~\ref{lem_C1} would lead to contradiction.

Let us prove that $r'_n \tend r_3$.
Otherwise, for a subsequence, we would have $r'_n \tend r' \in
[r_1,r_3[$. The sequence of holomorphic functions $z\mapsto
\phi_{\theta_n} (r'_n z)$ would converge uniformly on compact sets of
$\frac{r_3}{r'} \D$ to $\phi_\theta(r'z)$, which would imply uniform
convergence on $\U$ of all the derivatives. According to lemma~\ref{lem_C1},
the curves $\phi_{\theta_n}(r'_n\U)$ would be
$(K+\frac{1}{2})$-quasicircles for $n$ big enough, which leads to
contradiction.

As soon as $r'_n > r_1$, since for all $r\in [r_1,r'_n[$ the curve
$\phi_{\theta_n}(r\U)$ is a $(K+1)$-quasicircle, by continuity
the curve $\phi_{\theta_n}(r'_n\U)$ is a $(K+1)$-quasicircle.

Lemma~\ref{lem_geom_2} then implies that $t\in\T \mapsto
\phi_{\theta_n}(r'_n,t)-\phi_{\theta}(r_3,t)$ 
tends uniformly to $0$ when $n\tend +\infty$.
Since $r'_n \tend r_3$, $t\in\T \mapsto
\phi_{\theta}(r'_n,t)-\phi_{\theta}(r_3,t)$ 
also tends uniformly to $0$, which yields condition~\ref{item_inf}.
\qed\medskip

Having $r'=r_1=r_2$ would make the sequel simpler, and enable to fix
any value $<r(\theta)$ for $r(\theta')$ in theorem~\ref{thm_main}, but
it seems to require more than lemma~\ref{lem_abc}.

\

The following lemma recalls elementary properties of linearization of
Siegel disks.

\begin{lemma}\label{lem_limits}
  Let us note $(r,t)$ for $r \exp(i 2\pi t)$.
  Assume $r_n>0$, $r_n \tend r$, $\theta \in \R$ and $\theta_n \tend
  \theta$ is a sequence of real numbers such that $r(\theta_n)>r_n$, and the
  maps $t\in\T \mapsto \phi_{\theta_n}(r_n,t)$ form a
  Cauchy sequence for $\|\cdot\|_\infty$. Let $\psi$ be the limit.
  Then $r(\theta) \geq r$, and $(r,t) \mapsto \psi(t)$ extends continuously
  $\phi_{\theta}$ to the closure of $B(0,r)$
  (if $r(\theta) > r$, this just means $\phi_{\theta}(r,t) =
  \psi(t)$).
  Moreover, $\psi$ is injective. In the case $r(\theta) =r$, then
  $\partial \Delta(\theta) = \psi(\T)$.
\end{lemma}

Let us now state the main result.

\newcounter{savedcounter}
\begin{theorem}\label{thm_main}
  For all Bruno number $\theta$, for all $r<r(\theta)$
  and all $\epsilon>0$, there exists a Bruno number $\theta'$ such
  that
  \begin{enumerate}
    \item $|\theta'-\theta| < \epsilon$
    \item $r<r(\theta')<r+\epsilon$
    \item $\phi_{\theta'}$ has a continuous extension $\psi$ to $r(\theta')\U$
    \item $\big\|\psi - \phi_\theta\big\|_\infty <
    \epsilon$ on the circle $r(\theta')\U$ \label{item_norminf}
    \item the boundary $\psi(r(\theta')\U)$ of $\Delta(\theta')$ is a
    Jordan curve
    \item it does not contain the critical point \label{item_nococri}
    \setcounter{savedcounter}{\theenumi}
    \item it is not a quasicircle
  \end{enumerate}
\end{theorem}
\proof It is enough to prove the claim without point
\ref{item_nococri}, because it is implied by point \ref{item_norminf}
for $\epsilon$ small enough.

\

\noindent We are going to define by induction a sequence $\theta_n$ of
parameters, an increasing sequence $r'_n \geq r$, and angles $u_n$, $v_n$,
coordinates of points on the circle $\U$ via the map $t \mapsto
\exp(i2\pi t)$.

\

\noindent The induction hypothesis will be $H_n$:
\begin{itemize}
\item $r(\theta_n) >r'_n$
\item for all $k\leq n$,
$\on{pinch}(L, x_k, y_k) > k$
\end{itemize}
where
$L = \phi_{\theta_{n}} (r'_n \U)$,
$x_k = \phi_{\theta_{n}}(r'_n e^{i2\pi u_{k}})$ and
$y_k=\phi_{\theta_{n}}(r'_n e^{i2\pi v_{k}})$ (the second condition is
empty for $n=0$ and $1$).

\

\noindent Let $\theta_0 = \theta$, $r'_0=r$ and $u_0$, $v_0$ be any
distinct angles.

\

\noindent For $n\geq 1$, assume that $\theta_k$, $r'_k$, $u_k$, $v_k$ are
defined for $0\leq k < n$, and that $H_{n-1}$ holds.
There exists a $\eta>0$ such that for all continuous injective
$\psi : r'_{n-1}\U \to \C$, the condition $\|\psi -
\phi_{\theta_{n-1}}\|_\infty < \eta$ on $r'_{n-1}\U$ implies
that for all $k<n$, the pinching parameterized 
by angles $u_k$, $v_k$ remains $>k$ for the Jordan curve
$\psi(r'_{n-1}\U)$.
Let us note $(r,t)$ for $r \exp(i 2\pi t)$.
Let $r_1=r'_{n-1}$ and $r_2$ such that $r_1<r_2<r(\theta_{n-1})$,
close enough to $r_1$ so that
\[\forall r'\in ]r_1,r_2[,\quad \sup_{t\in\R} \big|
\phi_{\theta_{n-1}}(r',t) - \phi_{\theta_{n-1}}(r_1,t)
\big| < \max(\eta,\epsilon/2^n)/2.\]
Let $\theta_n$ and $r'_n$ be provided by lemma \ref{lem_maillon} such that
\begin{itemize}
  \item $|\theta_n-\theta_{n-1}|< \epsilon/2^n$
  \item $r_1 < r'_n <r_2$
  \item $r(\theta_n) > r'_n $
  \item $\big\|\phi_{\theta_n} - \phi_{\theta_{n-1}}\big\|_\infty <
  \max(\eta, \epsilon/2^n)/2$ on the circle $r'_n\U$
  \item the curve $\phi_{\theta_n} (r'_n\U)$ has a new pinching $>n$
\end{itemize}
We then define $u_n$ and $v_n$ as the angles parameterizing the new pinching.

\

\noindent Now that the sequences have been defined, let $\theta'$ be
the limit of the Cauchy sequence $\theta_n$, and $r'$ the limit of the
increasing sequence $r'_n$ (which is bounded from above by $4$).
Let us recall that for all $n$, $r(\theta_n)>r'_n$, and that
the sequence of maps $t\in\T \mapsto \phi_{\theta_n}(r'_n,t)$ is a
Cauchy sequence, whose limit we will call $\psi$. Thus we can apply
lemma \ref{lem_limits}: $\psi(r',t)$ continuously extends $\phi_{\theta'}$
to the closed ball $\overline{B}(0,r')$. Moreover, $\psi$ is injective,
thus for all $k \in \N$, the pair $(u_k,v_k)$ parameterizes a pinching
of the Jordan curve $L = \psi(\T)$, with $\on{pinch} \geq k$ by
continuity. Therefore $L$ is not a quasicircle, and $r(\theta')$
cannot be $>r'$ (otherwise, $L$ would be an analytic curve). Thus
$L=\partial \Delta(\theta')$.
\qed

\

\noindent\textsl{\large Variation}

\

A variation yields the next stronger theorem.

\

Let us call modulus of continuity any non decreasing positive function
$h$ defined on $[0,+\infty[$ and such that $h(\eta) \tend
0$ when $\eta \tend 0$. A function $f$ between compact metric spaces
is said to have $h$ as modulus of continuity if and only if
$d\big(f(x),f(y)\big) < h\big(d(x,y)\big)$ for all pairs $(x,y)$ with
$x\not=y$. We will say that $f$ 
is $h$-regular if there is $\lambda>0$ such that $f$ has $\lambda h$
as a modulus of continuity. For $g$, not being $h$-regular is
equivalent to: there exists sequences $x_n\not=y_n$ with
$d(x_n,y_n)\tend 0$ and $d\big(f(x),f(y)\big)/h\big(d(x,y)\big) \tend
+\infty$.

\begin{theorem}\label{thm_2}
  Let us make the same assumptions as in theorem~\ref{thm_main}.
  Let $h$ be any modulus of continuity. Then there exists
  a Bruno number $\theta'$ such that the same conclusions as in
  theorem~\ref{thm_main} hold, except for the following replacement:
  \begin{enumerate}
    \setcounter{enumi}{\value{savedcounter}}
    \item the map $\psi$ is not $h$-regular
  \end{enumerate}
\end{theorem}
If one takes (for instance) $h(\eta) = 1/|\log \eta|$, this implies
$\partial \Delta(\theta)$ is not a quasicircle, because the conformal
map of a quasidisk is always Hölder-continuous.

\

To prove theorem~\ref{thm_2}, we need to adapt lemma~\ref{lem_geom_2}:

\begin{lemma}\label{lem_geom_3}
Let $L_n$ be Jordan curves, $U_n$ be the bounded component of
$\C-L_n$. Assume that $0\in U_n$ and let $f_n : \D \to U_n$ be the
conformal isomorphism mapping $0$ to $0$ with real positive derivative.
Let $g_n$ be the continuation of $f_n$ to $\overline{\D}$ (exists
since $L_n$ is locally connected). Assume the retrictions of $g_n$ to
$\U$ have a common modulus of continuity, and 
that $U_n$ has Carathéodory limit equal to $\D$. Then 
the $g_n$ tend uniformly to identity.
\end{lemma}
\proof
According to Ascoli's theorem, the equicontinuous family
$g_n\big|_\U$ lies in a compact family of $C(\U)$ (the set of
continuous functions on $\U$ with the supremum norm). By the maximum
principle, for all $m,n\in\N$, 
the supremum of $|g_n-g_m|$ on $\overline{\D}$ is equal to its supremum on
$\U$. So, $g_n$ lies in a compact family of $C(\overline{\D})$.
The Carathéodory convergence of $U_n$ to $\D$ states that $g_n$ tends
to $\on{id}_\D$ uniformly on compact subsets of $\D$. So
$\on{id}_{\overline{\D}}$ is the only possible uniform limit on
$\overline{\D}$ of subsequences of $g_n$.
\qed

\

Lemma~\ref{lem_C1}

\begin{lemma}\label{lem_h}
  If $h$ is the modulus of continuity of a non constant function
  $f : \T \to \C$, then
  \[\inf_{\eta\in ]0,1]} \frac{h(\eta)}{\eta} > 0\]
  Therefore, if $f,f_n : \T \to \C$ are $C^1$ functions and such that 
  $f_n-f$ and its derivative uniformly tend to $0$,
  and $f$ has modulus of continuity $h$, then, for all $\epsilon>0$,
  $f_n$ has eventually modulus $(1+\epsilon) h$.
\end{lemma}

\

And lemma~\ref{lem_maillon}:

\begin{lemma}\label{lem_maillon_2}
  For all $\theta$ Bruno, for all $r_1<r_2<r(\theta)$, for all $K>1$ and for
  all $\epsilon>0$ there exists a Bruno number $\theta'$ and $r'>0$ such that
  \begin{enumerate}
   \item $|\theta'-\theta|<\epsilon$
   \item $r_1<r'<r_2$
   \item $r(\theta') > r'$
   \item \label{item_inf_2} $\| \phi_{\theta} - \phi_{\theta'}
   \|_\infty < \epsilon$ on
   the circle $r'\U$
   \item the restriction of $\phi_{\theta'}$ to $r'\U$ has not
   $h$ as a modulus of continuity
   \label{item_mc}
  \end{enumerate}
\end{lemma}
\proof
The proof is a straightforward adaptation of that of
lemma~\ref{lem_maillon_2}. However, we include it:\\
Let $r_3 = \frac{r_1+r_2}{2}$.
If, for some $r'\in ]r_1,r_3[$, $z\in \U \mapsto \phi_{\theta}(r' z)$
has not modulus 
of continuity $h$, then we are done with $\theta'= \theta$. 
Otherwise, let $\theta_n \tend \theta$ provided by lemma~\ref{lem_abc}
such that $r(\theta_n) \tend r_3$.
We claim that for all $n$ big enough, there exists an $r\in
[r_1,r(\theta_n)[$ such that 
$z\in \U \mapsto \phi_{\theta_n}(r z)$ has not modulus of continuity $2h$.
Otherwise, lemma~\ref{lem_geom_3} would imply that
$\partial \Delta(\theta_n)$ tends to $\phi_\theta(r_3 \U)$ and thus is
eventually contained in $\Delta(\theta)$, contradicting
lemma~\ref{lem_Herman} since $\theta_n \in \cal D_2$.

Let then $r'_n$ be the infimum of the set of $r \in [r_1,r(\theta_n)[$ such
that $\phi_{\theta_n} (r'_n z)$ has not modulus $2h$ on $\U$. 
Then according to lemma~\ref{lem_h}, $\phi_{\theta_n} (r z)$
has not modulus $\frac{3}{2} h$ on $\U$. 

Let us prove that $r'_n \tend r_3$.
Otherwise, for a subsequence, we would have $r'_n \tend r' \in
[r_1,r_3[$. The sequence of holomorphic functions $z\mapsto
\phi_{\theta_n} (r'_n z)$ would converge uniformly on compact sets of
$\frac{r_3}{r'} \D$ to $\phi_\theta(r'z)$, which would imply uniform
convergence on $\U$ of all the derivatives. Therefore, by lemma~\ref{lem_h},
$\phi_{\theta_n}(r'_n z)$ would eventually have modulus $\frac{3}{2} h$,
leading to contradiction.

As soon as $r'_n > r_1$, since for all $r\in [r_1,r'_n[$ the curve
$\phi_{\theta_n}(r z)$ has modulus $2h$ on $\U$, so does the function
$\phi_{\theta_n}(r'_n z)$.

Lemma~\ref{lem_geom_3} then implies that $t\in\T \mapsto
\phi_{\theta_n}(r'_n,t)-\phi_{\theta}(r_3,t)$ 
tends uniformly to $0$ when $n\tend +\infty$.
Since $r'_n \tend r_3$, $t\in\T \mapsto
\phi_{\theta}(r'_n,t)-\phi_{\theta}(r_3,t)$ 
also tends uniformly to $0$, which yields condition~\ref{item_inf_2}.
\qed

\

\part{Siegel disks with prescribed regularity}

\

We now inspire from the presentation in \cite{A} to give the following
theorem. Let $C^0$ be the space of holomorphic functions from $\D$ to
$\C$ having a continuous extension to $\overline{\D}$. This is a
Banach space for the supremum norm. Let $C^\omega$
be the space of functions from $\D$ to $\C$ having a holomorphic
extension to a neighborhood of $\D$. This is not a Fréchet space. This
is the union of spaces $C^\omega_\epsilon$ for $\epsilon >0$, where
$C^\omega_\epsilon$ is the set of holomorphic functions on
$(1+\epsilon)\D$. These spaces are endowed with the topology of uniform
convergence on compact sets. We do not put a topology on $C^\omega$.

\begin{theorem}\label{thm_control}
  Let us make the same assumptions as in theorem~\ref{thm_main}.
  Let $F$ be a Fréchet space such that
  \[C^\omega \subset F \subset_{\on{0}}
  C^0,\] where $\subset_{\on{0}}$ means a continuous injection.
  Assume that $K_n$ are compact subsets of $F$.
  Then there exists a Bruno number $\theta'$ such that the same
  conclusions as in theorem~\ref{thm_main} hold, except for the
  following replacement:
  \begin{enumerate}
    \setcounter{enumi}{\value{savedcounter}}
    \item the map $\psi$ belongs to $F$ but to no $K_n$
  \end{enumerate}
\end{theorem}

Note: a decreasing intersection of Fréchet spaces is a Fréchet
space. That is why there is no $F_n$ in the statement.

\

By the way, this proves that the set $\ds F - \bigcup K_n$ is not
empty ~! (This is of course very classical: Baire's theorem implies
it, since every compact subset of an infinite dimensional Fréchet
space has empty interior.) Even better~: it contains a univalent map.

\

Let's do the proof:

\begin{lemma}\label{lem_geom_4}
Let $F\subset_0 C^0$ be a Fréchet space
and $K$ be a compact subset of $F$.
Let $L_n$ be Jordan curves, $U_n$ be the bounded component of
$\C-L_n$. Assume that $0\in U_n$ and let $f_n : \D \to U_n$ be the
conformal isomorphism mapping $0$ to $0$ with real positive derivative.
Let $g_n$ be the continuation of $f_n$ to $\overline{\D}$ (exists
since $L_n$ is locally connected). Assume that $\forall n\in\N$, $g_n
\in K$, and that $U_n$ has Carathéodory limit equal to $\D$.
Then \par
a) $d_F(g_n,\on{id}_\D) \tend 0$
where $d_F$ is the distance
function of the Fréchet space $F$, 
\par b) $\|g_n-\on{id}_\D \|_\infty \tend 0$.
\end{lemma}
\proof
Part b) is a corollary of part a) and of the continuity of the injection
$F\subset C^0$.
Since $g_n$ lies in a compact set, it is enough to prove that all
convergent subsequences tend to $\on{id}_\D$. So we may assume
$d_F(g_n,h) \tend 0$ for some $h\in K$. The injection $F \subset C^0$
being continuous, $\|g_n-h\|_\infty \tend 0$. Carathéodory convergence
means that $g_n$ tends to $\on{id}_\D$ uniformly on compact sets. Thus 
$h=\on{id}_\D$.
\qed\medskip

\begin{lemma}\label{lem_abs}
 Assume $F$ is a Fréchet space such that 
 \[C^\omega \subset F \subset_0 C^0\]
 then, for all $\epsilon$, the injection $C^\omega_\epsilon \subset F$ is
 continuous.
\end{lemma}
\proof Corollary of the closed graph theorem, since the injection of
$C^\omega_\epsilon$ in $C^0$ is continuous.
\qed\medskip

\

\begin{lemma}\label{lem_B_n}
 There exists subsets $B_n$ of $C^\omega$ such that $C^\omega =
 \bigcup B_n$ and for all Fréchet space $F$ with $C^\omega \subset F
 \subset_0 C^0$, $B_n$ is compact in $F$.
\end{lemma}
\proof Let $B_n$ be the set of holomorphic functions on $\D$ which
have a holomorphic extension to $(1+\frac{1}{n+1}) \D$, that is
bounded by $n$. Each $B_n$ is compact in $C^\omega_{\frac{1}{n+1}}$
(Montel's theorem). According to lemma~\ref{lem_abs}, $B_n$ is also
compact in $F$.
\qed\medskip

\

We will consider the following property of a subset $A$ of $C^0$:
\[(\cal H) \qquad \forall \epsilon>0,\ A 
\text{ contains a neighborhood of } 0 \text{ in } C^\omega_\epsilon.\]

\begin{lemma}\label{lem_L}
 For all Fréchet space $F$ with $C^\omega \subset F \subset_0 C^0$,
 there exists a compact subset $L$ of $F$ with property $\cal H$.
\end{lemma}
\proof
 Take the same $B_n$ as in the proof of lemma~\ref{lem_B_n}.
 Since $B_n$ is compact in $F$, $\lambda B_n$ tends to $0$ when
 $\lambda \tend 0$. Choose $\lambda_n>0$ such that $\lambda_n B_n$ is
 included in the ball of $F$ of center $0$ and radius $1/(n+1)$. Then 
 the set $L= \bigcup \lambda_n B_n$, which has obviously
 property $\cal H$, is compact.
\qed\medskip

\

\begin{lemma}[perturbation]\label{lem_maillon_3}
  Assume that $F$ is a Fréchet space with $C^\omega \subset F
  \subset_0 C^0$, and that $K$ is a compact subset of $F$.
  For all $\theta$ Bruno, for all $r_1<r_2<r(\theta)$,
  and all $\epsilon>0$ there exists a Bruno number $\theta'$ and
  $r'>0$ such that 
  \begin{enumerate}
   \item $|\theta'-\theta|<\epsilon$
   \item $r_1<r'<r_2$
   \item $r(\theta') > r'$
   \item $d_F((\phi_{\theta} - \phi_{\theta'}) (r'z)) <
   \epsilon$\label{item_inf_bis}
   \item $\phi_{\theta'} (r'z)$ does not belong to $K$\label{item_five}
  \end{enumerate}
\end{lemma}
\proof
Let $r_3 = \frac{r_1+r_2}{2}$. 
If, for some $r'\in ]r_1,r_3[$, $\phi_{\theta}(r'z)$ verifies
(\ref{item_five}), then we are done with $\theta'= \theta$. 
\\
Otherwise, let $\theta_n \tend \theta$ provided by lemma~\ref{lem_abc}
such that $r(\theta_n) \tend r_3$. We may assume that $K$ has property
$\cal H$ by replacing $K$ by $K \cup L$ where $L$ is provided by
lemma~\ref{lem_L}. 
According to
lemma~\ref{lem_abs}, the injection $C^\omega_\epsilon \subset F$ is
continuous for all $\epsilon>0$.
On the other hand, there is some $\epsilon>0$ such that the function
$[r_1,r_3] \to C^\omega_\epsilon$ which maps $r$ to the function
$\phi_\theta(r z)$, is continuous. Therefore, its image $K_0$ is
compact in $C^\omega_\epsilon$, thus compact in $F$.
Let $K' = K + K_0$.
We claim that for all $n$ big enough, there exists an $r\in
[r_1,r(\theta_n)[$ such that 
$\phi_{\theta_n}(r z)$ does not belong to $K'$.
Otherwise, $\phi_{\theta_n}(r(\theta_n) z)$ would belong to $K'$ and
lemma~\ref{lem_geom_4}~b) would imply that 
$\|\phi_{\theta_n}(r(\theta_n) z) - \phi_\theta(r_3 z) \|_\infty \tend 0$
and thus eventually,
$\partial \Delta(\theta_n)$ would be contained in $\Delta(\theta)$,
contradicting lemma~\ref{lem_Herman} since $\theta_n \in \cal D_2$.

Let then $r'_n$ be the infimum of the set of $r \in [r_1,r(\theta_n)[$ such
that $\phi_{\theta_n} (r z)$ does not belong to $K'$. 

Let us prove that $r'_n \tend r_3$.
Otherwise, for a subsequence, we would have $r'_n \tend r' \in
[r_1,r_3[$. With the definition of $r'_n$, this would yield a sequence
$r''_n \tend r'$ with $r''_n>r'_n$ and
$\phi_{\theta_n}(r''_n z)$ does not belong to $K'$.
The sequence of holomorphic functions $z\mapsto
\phi_{\theta_n} (r''_n z)$ would converge uniformly on compact sets of
$\frac{r_3}{r'} \D$ to $\phi_\theta(r'z)$, which would imply uniform
convergence on $\U$ of all the derivatives. Because of property $\cal H$,
the function
$(\phi_{\theta_n}-\phi_\theta)(r''_n z)$ would
eventually belong to $K$, and
$\phi_{\theta_n}(r''_n z)$ would belong to $K+K_0$, that is $K'$, which
is a contradition.

As soon as $r'_n > r_1$, since for all $r\in [r_1,r'_n[$,
$\phi_{\theta_n}(r z)$ is in $K'$, by continuity
so is $\phi_{\theta_n}(r'_n z)$.

Lemma~\ref{lem_geom_4}~a) then implies that
$d_F(\phi_{\theta_n}(r'_n z)-\phi_{\theta}(r_3 z)) \tend 0$ when
$n\tend +\infty$. And $r'_n \tend r_3$ implies
$d_F(\phi_{\theta}(r_3 z)-\phi_{\theta}(r'_n z)) \tend 0$. This
gives~(\ref{item_inf_bis}).
\qed\medskip

\

\noindent\textsl{Proof of theorem~\ref{thm_control}:}

\

Let $B_n$ be provided by lemma~\ref{lem_B_n}, and $L$ by lemma~\ref{lem_L}.

\

We are going to define by induction a sequence $\theta_n$ of
parameters, an increasing sequence $r'_n \geq r$, and reals $\epsilon_n>0$.

\

\noindent The induction hypothesis will be $H_n$:
\begin{itemize}
\item $r(\theta_n) >r'_n$
\item for all $k\leq n$, the $F$-distance between 
$\phi_{\theta_n}(r'_n z)$ and the set $K_k\cup B_k \cup L$ is $>\epsilon_k$
\end{itemize}

\

\noindent Let $\theta_0 = \theta$, $r'_0=r$.

\

\noindent For $n\geq 1$, assume that $\theta_k$, $r_k$, $\epsilon_k$ are
defined for $0\leq k < n$, and that $H_{n-1}$ holds.
There exists a $\eta>0$ such that for all $f\in F$, the condition
$d_f(f,\phi_{\theta_{n-1}}(r'_{n-1}) ) < \eta$ implies
that for all $k<n$, the $F$-distance between $f$ and $K_k \cup B_k
\cup L$ remains
$>\epsilon_k$.
Let $r_1=r'_{n-1}$ and $r_2$ such that $r_1<r_2<r(\theta_{n-1})$,
close enough to $r_1$ so that
\[d_F(\phi_{\theta_n-1}(r' z), \phi_{\theta_n-1}(r_1 z)) < \max(\eta, \epsilon/2^n)/2\]
(possible since the injection of $C^\omega_{\epsilon'} \subset F$ is
continuous for all $\epsilon'>0$)
Let $\theta_n$ and $r'_n$ be provided by lemma \ref{lem_maillon_3} such that
\begin{itemize}
  \item $|\theta_n-\theta_{n-1}|< \epsilon/2^n$
  \item $r_1 < r'_n <r_2$
  \item $r(\theta_n) > r'_n $
  \item $d_F((\phi_{\theta_n} - \phi_{\theta_{n-1}})(r'_n z)) <
  \max(\eta, \epsilon/2^n)/2$
  \item $\phi_{\theta_n} (r'_n z)$ does not belong to $K_n \cup B_n \cup L$
\end{itemize}
We then define $\ds \epsilon_n = \frac{1}{2}d_F(K_n \cup B_n \cup L,\phi_{\theta_n} (r'_n z))$.

\

\noindent Now that the sequences have been defined, let $\theta'$ be
the limit of the Cauchy sequence $\theta_n$, and $r'$ the limit of the
increasing sequence $r'_n$ (which is bounded from above by $4$).
Let us recall that for all $n$, $r(\theta_n)>r'_n$, and that
the sequence of maps $\phi_{\theta_n}(r'_n z)$ (restricted to $\D$) is a
Cauchy sequence for $d_F$, thus converges in $F$ (that is where the
completeness of Fréchet spaces is used). Its limit is
$\psi(z)=\phi_{\theta'}(r' z) \in F$ (a priori restricted to $\D$).
Convergence in $F$ implies convergence in $C^0$, thus we can apply
lemma~\ref{lem_limits}. Also, $d_F( K_n \cup B_n \cup L , \psi ) \geq
\epsilon_n$, thus $\psi$ does not belong to any $K_n$ nor to any
$B_n$. Since $\bigcup B_n = C^\omega$, this implies $\psi$ does not
extend holomorphically to a neighborhood of $\D$, thus $r(\theta')=r'$.
\qed

\

\noindent\textsl{Examples:}

\

To obtain Siegel disks with smooth ($C^\infty$) boundaries, one takes
$F=C^\infty$ and $K_n = \emptyset$.

Let $B$ be a banach space (or a Fréchet space), and assume that 
\[C^\omega \subset B \subset_{\on{c}} F\]
where $\subset_{\on{c}}$ means a compact injection (the image of a bounded
set has compact closure). 
If we take $K_n$ to be the closure in $F$ of the 
the ball in $B$ of center $0$ and radius $n+1$,
we obtain Siegel disks whose conformal map $\phi_{\theta'} (r'z)$ 
belongs to $F$ but not to $B$. For instance
\begin{itemize}
 \item $F=C^0$, $B$ is the set of functions whose restriction to
 $\U$ is $h$-regular: this
 reproves theorem~\ref{thm_2}
 \item $F=C^n$, $B=C^{n+1}$ where $n\in\N$
 \item $\ds F=\bigcap_{\alpha \in [0,\beta[} C^\alpha$, $B=C^\beta$, where
 $\beta > 0$ is a real number
 \item $\ds F=\bigcap_{\alpha \in [0,n[} C^\alpha$, $B=C^{(n-1)+\on{Lip}}$, where
 $n > 0$ is an integer
\end{itemize}
In the first two examples, the Fréchet space $F$ happens to be a
Banach space.

Now we can take a countable collection of
Banach spaces
$B_n$ such that $C^\omega \subset B_n \subset_{\on{c}} F$. We obtain
Siegel disks whose conformal map $\phi_{\theta'} (r'z)$ 
belongs to $F$ but to no $B_n$. For instance
\begin{itemize}
 \item $\ds F=C^\alpha$, $B_n=C^{\alpha+1/n}$, where $\alpha\geq 0$
\end{itemize}

\remark Since the inclusion of $C^n$ in $C^{(n-1)+\on{Lip}}$ is
not compact (it is an isometry), we may
wonder if there exists Siegel disks whose boundaries are
$C^{(n-1)+\on{Lip}}$ but not $C^n$. For $n=1$, the contrary would mean 
that if the boundary of a (fixed quadratic) Siegel disk is $\on{Lip}$,
then it is $C^1$.


\begin{thebibliography}{McM2}

\bibitem[A] {A} {\sc A.\ Avila}, {\em Smooth Siegel disks via
semicontinuity: a remark on a proof of Buff and Cheritat},
\texttt{math.DS/0305272}

\bibitem[BC] {BC1} {\sc X.\ Buff, A.\ Ch\'eritat}, {\em Quadratic
Siegel disks with smooth boundaries. Part~I}, submited.

\bibitem[H]{He} {\sc M.R.\ Herman}, \emph{Are there critical
points on the boundaries of singular domains~?}
Comm.\ Math.\ Phys.\ \textbf{99}, 593--612 (1985).

\bibitem[P]{p} {\sc C.\ Pommerenke}, {\em Boundary Behavior of Conformal Maps},
Grundlehren der mathematischen Wissenschaften~299, Springer-Verlag.

\bibitem[PZ] {PZ} {\sc  C.L. Petersen} $\&$ {\sc S. Zakeri}, {\em On the Julia
Set of a Typical Quadratic Polynomial with a Siegel disk}, Preprint,
Institute for Mathematical sciences, SUNY at Stony Brook,  (2000).

\end{thebibliography}
\end{document}